\documentclass[11pt]{amsart}
\usepackage{amsmath}
\usepackage{latexsym}
\usepackage{amssymb,amsthm,amsfonts}
\usepackage{graphicx}

\newtheorem{defi}{Definition}[section]
\newtheorem{theorem}[defi]{Theorem}
\newtheorem{lemma}[defi]{Lemma}

\newtheorem{proposition}[defi]{Proposition}

\newenvironment{Proof}{\noindent{\bf Proof}.}{\hfill$\Box$\\[2mm]}

\newcommand{\beq}{\begin{equation}}
\newcommand{\eeq}{\end{equation}}
\newcommand{\bdoc}{\begin{document}}
\newcommand{\edoc}{\end{document}}
\newcommand{\be}{\begin{enumerate}}
\newcommand{\ee}{\end{enumerate}}
\newcommand{\bdescr}{\begin{description}}
\newcommand{\edescr}{\end{description}}
\newcommand{\ba}{\begin{array}}
\newcommand{\ea}{\end{array}}

\def\R{\mathbb{R}}
\def\e{\varepsilon}
\def\dis{\displaystyle}
\def\intr{\int_{\R^3}}

\title[Cluster solutions for the Schr{\"o}dinger-Poisson-Slater problem]{Cluster solutions for the Schr{\"o}dinger-Poisson-Slater problem around a local minimum of the potential}
\author{David Ruiz and Giusi Vaira}

\address{Departamento de An\'alisis Matem\'atico, University of Granada, 18071 Granada (Spain) and SISSA, via Beirut 2-4, 34014 Trieste (Italy)}

\thanks{D.R has been supported
by the Spanish Ministry of Science and Technology under Grant
MTM2005-01331 and by J. Andaluc\'{\i}a (FQM 116).}

\email{daruiz@ugr.es, vaira@sissa.it}


\subjclass[2000]{}

\date{}

\begin{document}

\maketitle

\begin{abstract}

In this paper we consider the system in $\R^3$
\begin{equation}\label{problemadipartenza0}
\left\{
  \begin{array}{l}
    -\e^2\Delta u+V(x)u+\phi(x)u=u^{p},

    \\

    -\Delta \phi=u^{2},
  \end{array}
 \right.
\end{equation}
for $p\in (1, 5)$. We prove the existence of multi-bump solutions
whose bumps concentrate around a local minimum of the potential
$V(x)$. We point out that such solutions do not exist in the
framework of the usual Nonlinear Schr{\"o}dinger Equation.

\end{abstract}

\section{Introduction and main results}
\noindent Recently, many papers have studied different versions of the
Schr{\"o}dinger-Poisson-$X^{\alpha}$ problem:

\begin{equation}\label{eq11} - \frac{\hbar^2}{2m}\Delta u + V(x) u + \left ( u^2 \star \frac{1}{4\pi|x|} \right )
u= |u|^{p-1}u, \ x \in \R^3, \end{equation} where $V(x)$ is an
external potential and $p \in (1,5)$. The interest on this problem
stems from the Slater approximation of the exchange term in the
Hartree-Fock model, see \cite{slater}. In this framework, $p=5/3$;
however, other exponents have been used in different
approximations, which have been referred to as $X^{\alpha}$ type
approximations, see \cite{mauser}. From another point of view,
this equation has been proposed in \cite{bfortunato} under the
name of Schr{\"o}dinger-Maxwell equation. For more information on the
relevance of this model and its deduction, we refer to
\cite{bfortunato, boka, bls, cornean, mauser}.

From the mathematical point of view, problem \eqref{eq11} presents
an interaction between two different kind of nonlinear terms: a
repulsive nonlocal term and an attractive local term. This, and
related problems, have been much studied recently by using
variational methods, see \cite{a-ruiz, azz, mugnai, mugnai2,
kikuchi, kikutesis, pisani, JFA, zhou, zhao}.

If we define $\phi_u= u^2 \star \frac{1}{4\pi|x|}$ and $\e^2=
\frac{\hbar^2}{2m}$, the equation \eqref{eq11} can be rewritten as
a system in the form:

\begin{equation}\label{problemadipartenza}
\left\{
  \begin{array}{l}
    -\e^2\Delta u+V(x)u+ \phi(x)u=|u|^{p-1}u,
    \\
    -\Delta \phi=u^{2}.
  \end{array}
 \right.
\end{equation}
In this paper we are concerned with the semiclassical limit for
the system (\ref{problemadipartenza}), namely the problem of
finding non trivial solutions $(u, \phi)\in H^{1}(\mathbb
R^3)\times D^{1, 2}(\mathbb R^3)$ and studying their asymptotic
behavior as $\e \rightarrow 0$. Such solutions are usually
referred to as {\it semiclassical states}. 

A large number of papers deals with the study of semiclassical
states for the following nonlinear Schr\"{o}dinger equation
\begin{equation}\label{Schrodinger}
-\e^2\Delta u+V(x)u=|u|^{p-1}u, \qquad x\in \mathbb R^3.
\end{equation}
For the problem (\ref{Schrodinger}) spike solutions are found
around the critical points of the potential $V$, see for instance
\cite{AmbrBadiCing, Li}. These are solutions that concentrate (as
$\e \to 0$) around a unique point, and tend to zero outside of
this point. For instance in \cite{Li} Yanyan Li proved the
existence of positive solutions concentrating near $C^1$ stable
critical points of $V$. Moreover, Li proves also the existence of
multi-bump solutions, namely, solutions concentrating around
different critical points of $V$. Other results in this direction
were given in \cite{delpino, gui}. However, in the previous papers
the bumps are well separated and so the interactions among the
different bumps are neglected.

In \cite{KangWei} the authors prove the existence of multi-bump
solutions for (\ref{Schrodinger}) whose bumps tend to a point of
local maximum of $V$. Here the interactions among the bumps do
play a role. In a certain sense, each bump has an attractive
effect on the other bumps, whereas the potential has a repulsive
effect (around its local maximum). The multi-bump solution exists
due to a balance between the two effects. The authors also show
that multi-bump solutions do not exist around nondegenerate local
minima. In this case, both effects
would be attractive and no balance could be possible.


With respect to (\ref{problemadipartenza}), the existence of
single-bump solutions near critical points of $V$ has been
recently proved, see \cite{IanniVaira}). Other concentration
phenomena have been proved for this system even with the absence
of the potential, see \cite{DaprileWei, DaprileWei2}.

In this paper we prove the existence of positive solutions with
$K$ intera-\ cting bumps around local minima of the potential $V$.
These solutions appear because of the effect of the Poisson term
in our equation. Indeed, the Poisson term implies a repulsive
effect among the bumps which balance the attractive effect of the
potential $V$.

We assume that:

\begin{itemize}
\item[(V1)] $V$ has a local strict minimum point in $P_0$, namely
there exists a bounded open set $\mathcal{U}$ such that $P_0\in
\mathcal{U}$ and $$V(P_0)=\min_{x\in
\bar{\mathcal{U}}}V(x)<V(P),\qquad \forall \,\, P\in
\mathcal{U}\setminus \{P_0\}$$ Up to a translation and dilatation,
we can assume $P_0=0$, $V(0)=1$. \item[(V2)] $V(x)= 1 +
|g(x)|^{\alpha}$ for any $x \in \mathcal{U}$, where $g:
\mathcal{U} \to \R$ is a $C^{2,1}$ function and $\alpha>2$.
\end{itemize}

In particular, there holds:

\begin{itemize} \item[(V2')] $V(x)\leq 1+C |x|^{\alpha}$ for $x\in
\mathcal{U}$ and some $C>0$. \end{itemize}

Observe that under the above conditions the local minimum must be
degenerate. We point out that conditions (V1)-(V2') are sufficient
for most of our arguments. We need condition (V2) for technical
reasons, to be able to rule out possible undesired oscillations of
the derivatives of $V$ near $0$.

Let us denote by $U$ the unique positive radial solution in
$H^1(\R^3)$ of the problem (see \cite{Kwong}):
\begin{equation} \label{P1} -\Delta U + U = U^p. \end{equation}

Our main result is the following.

\begin{theorem}\label{minimo}
Assume that $V$ satisfies (V1) and (V2) and suppose $p\in (1, 5)$.
Then for any positive integer $K \in \mathbb Z$, there exists
$\e_K>0$ such that for any $\e<\e_K$ there exists a positive
solution $u_{\e}$ of (\ref{problemadipartenza}) with $K$ bumps
converging to $0$. More specifically, there exists $Q_1^{\e},
\dots Q_k^{\e} \in \R^3$ such that:
\begin{enumerate}
\item $Q_i^{\e} \to 0$, $\e^{-1}|Q_i^{\e}| \to +\infty$  as $\e
\to 0$.

\item Defining $\tilde{u}_{\e}(x)=u_{\e}(\e x)$, we have that
$\tilde{u}_{\e}(x) = \sum_{i=1}^K U(x-\e^{-1} Q_i^{\e}) + o(1)$,
as $\e \to 0$.

\end{enumerate}

\end{theorem}
The proof uses a singular perturbation method, based on a
Lyapunov-Schmidt reduction. We point out that the distance between
the bumps is different from that of the multi-bump
solutions of \cite{KangWei}, and this is caused because the
different balance involving the Poisson term.

The paper is organized as follows. Section \ref{Preliminari} is
devoted to some notations and to the variational setting of the
problem. In Section \ref{casominimo} we introduce the
Lyapunov-Schmidt reduction and solve the auxiliary equation.
Finally, in Section \ref{sezioneridotto} the reduced functional is
studied, solving the bifurcation equation. This completes the
proof of Theorem \ref{minimo}.

\section{\bf Preliminaries}\label{Preliminari}
\noindent As mentioned in the introduction, we denote by $U$ the unique
positive radial solution in $H^1(\R^3)$ of the problem
$$ -\Delta U + U = U^p. $$
This solution satisfies the following decay property (see
\cite{Kwong}):
\begin{displaymath}
\lim_{r\rightarrow +\infty} U(r) r e^{r} =C>0 , \,\ \,\ \,\
\lim_{r\rightarrow +\infty}\frac{U'(r)}{U(r)}=-1, \,\ \,\ \,\
r=|x|.
\end{displaymath}
for some constant $C$.

The function $U$ is a critical point of the $C^2$ functional
$I_{0}:H^{1}({\mathbb R}^{3}) \rightarrow {\mathbb R}$ defined as
\begin{equation}\label{I1}
I_{0}(u)=\frac{1}{2}\|u\|^2
                -\frac{1}{p+1}\int_{{\mathbb R}^{3}}|u|^{p+1}\, dx,
\end{equation}
where $\| \cdot \|$ denotes the usual norm in $H^1(\R^3)$.
Furthermore the solution $U$ is nondegenerate (up to
translations). More specifically, there holds:
\begin{lemma}\label{nondeg}
   {\it  Define the operator $Q:H^1(\mathbb R^3)\rightarrow \mathbb R$ as
    $$Q[\nu]:=I_{0}''(U)[\nu, \nu]=\int_{\mathbb R^3}\Big[|\nabla\nu|^2+\nu^2-pU^{p-1}\nu^2 \Big]\, dx.$$
    We denote $U_k=\frac{\partial U}{\partial x_k}$. Then there hold:
\begin{itemize}
    \item $Q[U]=(1-p)\|U\|^2<0$.
    \item $Q[\frac{\partial U}{\partial x_j}]=0,\,\,j=1, 2, 3.$
    \item $Q[\nu]\geq C\|\nu\|^2$ for all $\nu\bot U, \nu\bot \frac{\partial U}{\partial x_j}$, $j=1, 2, 3.$
\end{itemize}}
    \end{lemma}
For a proof see for instance \cite[Lemma 8.6]{AmbrMalc}.

It is convenient to make the change of variable $x \mapsto \e x$
and so we arrive to the problem:

\begin{equation}\label{problemafinale}
-\Delta u+V(\e x)u+\e^{2}\phi_{u}u=u^{p},\qquad u \in H^1(\mathbb
R^3),\qquad u>0.
\end{equation}

Here $\phi_u \in D^{1,2}(\R^3)$, and

$$ \int_{\R^3} |\nabla \phi_u|^2 \, dx= \int_{\R^3} \phi_u u^2\,dx=\int_{\R^3} \int_{\R^3}
\frac{u^2(x)u^2(y)}{4\pi|x-y|}\, \, dx \, dy.$$

In general, given $f \in L^{6/5}$, the solution of the problem
$-\Delta \phi = f$ belongs to $D^{1,2}(\R^3)$ and:
$$ \int_{\R^3} \nabla \phi \cdot \nabla \psi = \int_{\R^3} f \psi \leq \|\psi\|_{L^6} \| f\|_{L^{6/5}}
\leq C \| \psi \|_{D^{1,2}} \|f\|_{L^{6/5}}.$$ Therefore, $\| \phi
\|_{D^{1,2}} \leq C \| f\|_{6/5}$.

\noindent Moreover, it is well-known (see \cite{bfortunato}, for
example) that the solutions of (\ref{problemafinale}) correspond
to positive critical points of the $C^{2}$ functional
$I_{\e}:H^1({\mathbb R}^3) \rightarrow {\mathbb R}$,
\begin{equation}\label{funzionale}
I_{\e}(u)=\frac{1}{2}\int_{{\mathbb R}^{3}}\left[|\nabla
u|^{2}+V(\e x)u^{2}\right]\, dx +\frac{\e^{2}}{4}\int_{{\mathbb
R}^{3}}\phi_{u}(x)u^{2}\, dx
                -\frac{1}{p+1}\int_{{\mathbb R}^{3}}|u|^{p+1}\, dx.
\end{equation}

Finally, let us compute the derivatives of $V$. By using (V2)
\begin{equation} \label{second} \begin{array}{c} V_{x_i}(x)=
\alpha |g(x)|^{\alpha -2} g(x) g_{x_i}(x), \\ \\ V_{x_i \,
x_j}(x)= \alpha (\alpha-1)|g(x)|^{\alpha -2} g_{x_i}(x) g_{x_j}(x)
+ \alpha |g(x)|^{\alpha -2} g(x) g_{x_i \, x_j}(x).
\end{array}\end{equation}

In particular $V \in C^{2,\, \gamma}(\mathcal{U})$, where $\gamma
= \min \{1, \alpha-2\}$.

\section{The Lyapunov-Schmidt reduction. The auxiliary equation}\label{casominimo}

In this section we begin the Lyapunov-Schmidt for the proof of
Theorem \ref{minimo}. This will be made around an appropriate set
of ``approximating solutions". For any $K\in \mathbb N$, we define
$$\Lambda_{\e}=\left\{{\bf P}\in \mathbb R^{3K}: |P_i-P_j|\geq \e^{\frac{2-\alpha}{\alpha+1}+\delta}, \,\,
i\neq j,\,\ V(\e P_i)\leq
1+\e^{\frac{3\alpha}{\alpha+1}-\delta},\,\, \e P_i \in
\mathcal{U}\right\}$$ where $\delta>0$ is chosen small enough so
that $\frac{3\alpha}{\alpha+1}-\delta > 2$ (this is possible since
$\alpha >2$). Observe that $\frac{2-\alpha}{\alpha+1}+\delta < 0$
and $\Lambda_{\e}$ is not empty for $\e$ small enough.

Fix ${\bf P}=(P_1, ... , P_K) \in \Lambda_{\e}$. Setting
$z_{P_i}(x)=U(x-P_i)$, we define the manifold of ``approximate
solutions":
$$\mathcal{Z}=\left\{ z_{\bf{P}}(x) = \sum_{i=1}^{K}z_{P_i}(x) \, :
\ \, {\bf P} \in \Lambda_{\e}\right\}.$$

This section is devoted to the proof of the next result:

\begin{proposition} \label{prop} Assume that $V$ satisfies (V1) and (V2) and suppose $p\in (1, 5)$.
Then for any positive integer $K \in \mathbb Z$, there exists
$\e_K>0$ such that for any $\e<\e_K$ there exists a positive
solution $u_{\e}$ of \eqref{problemafinale}, and $z_{\e} \in
\mathcal{Z}$ such that $\|u_{\e} - z_{\e}\| = O(\e^2)$.
\end{proposition}

It is easy to check that Proposition \ref{prop} implies Theorem
\ref{minimo}.

The proof uses a Lyapunov-Schmidt reduction. For every $z \in
\mathcal{Z}$, we define
$W=W_{z,\e}=\left(T_{z}\mathcal{Z}\right)^{\bot}$ and $P:
H^1(\mathbb R^3) \rightarrow W$ the orthogonal projection onto
$W$. Our approach is to find a pair $z\in \mathcal{Z}$, $w \in W$,
$\| w \| = O(\e^2)$, such that $I_{\e}'(z+w) = 0$. Equivalently:
\begin{equation} \label{lya}
\left\{
  \begin{array}{lr}
    \mbox{a})\,\, PI_{\e}'(z+w)=0,
    \\
    \\
    \mbox{b})\,\, (\mathcal{I}-P)I_{\e}'(z+w)=0.
  \end{array}
 \right.
\end{equation}
The first equation above is called auxiliary equation, and the
second one receives the name of bifurcation equation.

\medskip

Our intention now is to find a solution $w \in W$ of the auxiliary
equation for any $z \in \mathcal{Z}$. We begin with some
estimates:

\begin{proposition}\label{soluzioniapprossimanti}
There exists $C=C(K)>0$ such that for all $\e>0$ small and any
${\bf P} \in \Lambda_{\e}$, we have
\begin{equation}\label{errore}
\|I'_{\e}(z_{{\bf P}})\| \leq C\e^2.
\end{equation}
\end{proposition}
\begin{Proof}
Taking into account that $z_{P_i}$ are solutions of (\ref{P1}), we
have:
\begin{eqnarray*}
I_{\e}'(z_{{\bf P}})[v]= \underbrace{\int_{\mathbb R^3}[V(\e
x)-1]z_{{\bf P}} v\, dx}_{(I)}+  \e^2\underbrace{\int_{\mathbb
R^3}\phi_{z_{{\bf P}}}z_{{\bf P}}v\,
dx}_{(II)}-\underbrace{\int_{\mathbb R^3}\left[|z_{{\bf P}}|^p-
\sum_{i=1}^K z_{P_i}^p\right]v\, dx}_{(III)}
\end{eqnarray*}
Let us evaluate separately the various terms. The second term can
be easily estimated (see Section \ref{Preliminari}):
\begin{eqnarray}\label{(II)}
(II) & \leq & \|z_{{\bf P}}\|^3\cdot \|v\| \leq C \, K \, \| v \|.
\end{eqnarray}

For (I), it suffices to estimate
\begin{eqnarray*}
 \int_{\R^3} [V(\e x)-1]z_{P_i}v \, dx \leq
\underbrace{\int_{\R^3}[V(\e x)-V(\e P_i)]z_{P_i}v\, dx}_{(A)}+
\underbrace{\int_{\R^3} [V(\e P_i)-1]z_{P_i}v\, dx}_{(B)}.
\end{eqnarray*}
By the definition of $\Lambda_{\e}$, we get that $(B)=o(\e^2)$.
Let us estimate $(A)$ by splitting the integral in two parts:

\begin{eqnarray*} \int_{\R^3}[V(\e x)-V(\e P_i)]z_{P_i}v\, dx = &
\dis \int_{|x-P_i|>\e^{-1}}[V(\e x)-V(\e P_i)]z_{P_i} v\, dx + \\
& \dis \int_{|x-P_i|<\e^{-1}}[V(\e x)-V(\e P_i)]z_{P_i} v\,
dx.\end{eqnarray*} Since $V$ is bounded in $L^{\infty}$, we use
H\"{o}lder estimate and the change $y=x-P_i$, to conclude
$$\int_{|x-P_i|>\e^{-1}}[V(\e x)-V(\e
P_i)]z_{P_i}v\, dx \leq C \left (\int_{|y|>\e^{-1}}U^2(y)\, dy
\right )^{1/2} \|v\|_{L^2} = o(\e^M) \|v\|_{L^2}$$ for any $M>0$,
thanks to the exponential decay of $U$.

Observe that if $|x-P_i|<\e^{-1}$, $\e P_i$ belongs to
$\mathcal{U}$ and $d(\e x, \mathcal{U}) \leq 1$. We use a Taylor
expansion:

\begin{equation} \label{ekeland} \begin{array}{cc} \dis \int_{|x-P_i|< \e^{-1}}|V(\e x)-V(\e P_i)| z_{P_i} |v| \, dx
\leq  \\ \\ \dis \int_{\R^3} \left (\e |\nabla V(\e P_i)| \,
|x-P_i| + C \e^2 |x-P_i|^2 \right )z_{P_i} |v| \, dx.
\end{array}
\end{equation}

Again by the exponential decay of $U$, $\| \, |x-P_i|^m z_{P_i}
\|_{L^2}$ is uniformly bounded for any $m >0$. So it suffices to
estimate $|\nabla V(\e P_i)|$.

Recall that $\e P_i \in \Lambda_{\e}$, and so $V(\e P_i)= 1+ |g(\e
P_i)|^{\alpha} \leq 1+ \e^{ \frac{3\alpha}{\alpha+1}-\delta}$. By
\eqref{second},
$$|V_{x_i}(x)| \leq C |g(x)|^{\alpha-1} \leq C \e^{ \left (\frac{3\alpha}{\alpha+1}-\delta \right )
\frac{\alpha-1}{\alpha}}.$$

Observe that $\frac{3\alpha}{\alpha+1}-\delta > 2 >
\frac{\alpha}{\alpha-1}$. Therefore, $\nabla V(\e P_i) = o(\e)$.

\medskip Finally we consider (III). These estimates have been done in
\cite{KangWei}; we sketch here the proof for the sake of
completeness. Let us define $\rho_{\e} =
\e^{\frac{2-\alpha}{\alpha+1}+\delta}$ and divide $\R^3$ in $K+1$
regions: $$\Omega_i=\{ x \in \R^3: 2 |x-P_i| \leq \rho_{\e}\} \ \
\mbox{ for } i=1 \dots K,\ \ \Omega_0 = \R^3 \setminus \left (
\cup_{i=1}^K \Omega_i \right ).$$

We now use the $C^{1,\sigma}$ regularity of the function
$f(u)=u^p$, where $\sigma = \min \{1, p-1\}$:

$$ \begin{array}{l} \dis \int_{\Omega_j}\left |\left ( \sum_{i=1}^K
z_{P_i} \right )^p- z_{P_j}^p - \sum_{i \neq j} z_{P_i}^p \right|
\, |v| \, dx\leq \\ \dis \int_{\Omega_j }\left[p z_{P_j}^{p-1} \left
(\sum_{i \neq j} z_{P_i} \right ) + C \left (\sum_{i \neq j}
z_{P_i} \right )^{1+\sigma} + \sum_{i \neq j} z_{P_i}^p
 \right]\, |v|\, dx \leq \\ \dis C \int_{\Omega_j } \left (\sum_{i \neq j}
z_{P_i} \right ) |v|\, dx. \end{array} $$ The last inequality is due to
the fact that in $\Omega_j$, $z_{P_i}\leq 1$. Indeed, defining
$\rho_{\e} = \e^{\frac{2-\alpha}{\alpha+1}+\delta}$ and using the
exponential decay of $U$, we have
$$ \int_{\Omega_j} z_{P_i}^2(x)\, dx \leq  \int_{2|x|> \rho_{\e}} U^2(y)\, dy \leq
C\int_{2r> \rho_{\e}} e^{-2 y}\, dr = C e^{-\rho_{\e}}.$$

On the other hand,

$$ \int_{\Omega_0} \left |\left ( \sum_{i=1}^K
z_{P_i} \right )^p- \sum_{i =1}^K z_{P_i}^p \right| \, |v| \leq C
\int_{\Omega_0} \sum_{i =1}^K z_{P_i}^p |v|,$$
$$ \int_{\Omega_0} z_{P_i}^{2p}(x) \, dx \leq  \int_{2|x|> \rho_{\e}} U^{2p}(y)\, dy \leq C e^{-p \rho_{\e}}. $$

This concludes the estimate (III).

\end{Proof}

Now we are concerned with the invertibility of $I''_{\e}(z_{\bf
P})$ on $W=(T_{z_{\bf P}}(\mathcal{Z}))^{\bot}$. First we observe
that $T_{z_{\bf P}}\mathcal{Z}$ is spanned by the functions
$\dot{z}_{i, j}:=\displaystyle\frac{\partial U}{\partial
x_j}(x-P_i)$, with $i=1,...,K$ and $j=1, 2, 3$. Recall that $P$
denotes the orthogonal projection onto $W$; me decompose: $W=A
\oplus B$ where
$$A = \langle \displaystyle \{ Pz_{P_i}\}_{i=1 \dots K} \rangle \
\mbox{and}\,\, B=\left(\displaystyle A \oplus T_{z_{\bf
P}}\mathcal{Z}\right)^{\bot}$$

\begin{proposition} \label{prop2} For $\e$ small and any ${\bf P} \in \Lambda_{\e}$, $PI_{\e}''(z_{\bf P}): W \to W$ is invertible and
$\|[PI_{\e}''(z_{\bf P})]^{-1}\| \leq \bar{C}$. \end{proposition}

The above result follows directly from the following lemma (see
\cite{AmbrMalc}):

\begin{lemma}\label{invertibilita}
For all $\e>0$ sufficiently small there exist two positive
constants $C_1, C_2$ such that
\begin{itemize}
\item[(a)] $I_{\e}''(z_{\bf P})[u, u] \leq -C_1 \| u \|^2$, for
all $u\in A$; \item[(b)] $I_{\e}''(z_{\bf P})[u, u] \geq C_2
\|u\|^2$, for all $u\in B$.
\end{itemize}
\end{lemma}
\begin{Proof}
Let be  $u\in A$. Then $$u=\sum_{i=1}^K\lambda_i Pz_{P_i},\qquad
\lambda_i\in \mathbb R,\quad i=1,..., K.$$ For $i=1,...,K$, $Pz_{P_i}$ are orthogonal to $T_{z_{\bf P}}(\mathcal{Z})$. Hence we can write
$$Pz_{P_i}=z_{P_i}-\psi_{i}, \qquad i=1,...,K$$ where $\psi_i$ are given by
$$\psi_{i}=\sum_{\begin{array}{c} l,\ j\\ l\neq i\end{array}}
\left(z_{P_i}, \dot{z}_{l, j}\right)\frac{\dot{z}_{l,
j}}{\left|\left|\dot{z}_{l, j}\right|\right|^2}.$$ The functions
$\dot{z}_{l, j}$ satisfy $-\Delta \dot{z}_{l, j}+\dot{z}_{l,
j}=pz_{P_l}^{p-1}\dot{z}_{l, j}$.\\ Since for $i\neq l$,
$|P_i-P_l|\rightarrow +\infty$ as $\e\rightarrow 0$, after an
integration by parts, we get $(z_{P_i},\dot{z}_{l, j})=o(1)$ as
$\e\rightarrow 0$. This implies $\|\psi_i\|=o(1)$ as
$\e\rightarrow 0$ for $i=1, ... , K$.

We now apply the bilinear form given by $I_{\e}''(z_{\bf P})$ to
obtain
\begin{eqnarray*}
I_{\e}''(z_{\bf P})[u, u]&=& \underbrace{I_{\e}''(z_{\bf P})\left[\sum_{i=1}^K\lambda_i z_{P_i}, \sum_{i=1}^K\lambda_iz_{P_i}\right]}_{(I)}+\underbrace{I''_{\e}(z_{\bf P})\left[\sum_{i=1}^K\lambda_i \psi_i, \sum_{i=1}^K\lambda_i\psi_i\right]}_{(II)}\\\\
&&+2\underbrace{I_{\e}''(z_{\bf P})\left[\sum_{i=1}^K\lambda_i
z_{P_i}\sum_{i=1}^K\lambda_i \psi_i\right]}_{(III)}.
\end{eqnarray*}
We observe that $I_{\e}''(z_{\bf P})$ maps bounded sets
onto bounded sets, then since $z_{\bf P}$ is bounded
\begin{eqnarray*}
(II)\leq \|I_{\e}''(z_{\bf P})\|\sum_{i=1}^K\lambda_i^2\|\psi_i\|^2\leq C\sum_{i=1}^K\lambda_i^2\|\psi_i\|^2=o(1).
\end{eqnarray*}
In the same way we obtain
\begin{eqnarray*}
(II)\leq \|I_{\e}''(z_{\bf P})\|\sum_{i=1}^K\lambda_i^2\|\psi_i\|^2\leq C\sum_{i=1}^K\lambda_i^2\|\psi_i\|=o(1).
\end{eqnarray*}
Furthermore, by making simple computations one finds
\begin{eqnarray*}
(I)&=&\sum_{i=1}^K \lambda_i^2\left(\int_{\mathbb R^3}[|\nabla z_{P_i}|^2+z_{P_i}^2-p z_{P_i}^{p+1}]\, dx\right)+\sum_{i=1}^K \lambda_i^2\underbrace{\left(\int_{\mathbb R^3}[V(\e x)-1]z_{P_i}^2\, dx\right)}_{(A)}\\\\
&&+2\sum_{i\neq j}\lambda_i\lambda_j\underbrace{\left(\int_{\mathbb R^3}[\nabla z_{P_i}\nabla z_{P_j}+V(\e x)z_{P_i}z_{P_j}]\, dx\right)}_{(B)}\\\\
&&+\underbrace{\e^2\int_{\mathbb R^3}\phi_{z_{\bf P}}\left(\sum_{i=1}^{K}\lambda_iz_{P_i}\right)^2\, dx}_{(C)}+
\underbrace{2\e^2\int_{\mathbb R^3} \widetilde{\phi} \cdot z_{\bf P}\left(\sum_{i=1}^K \lambda_i z_{P_i}\right)dx}_{(D)}\\\\
&&-\underbrace{p \int_{\mathbb R^3}\left[\left|\sum_{i=1}^K z_{P_i}\right|^{p-1}\left(\sum_{i=1}^K\lambda_iz_{P_i}\right)^2-\sum_{i=1}^K\lambda_i^2z_{P_i}^{p+1}\right]\, dx}_{(E)}\\
\end{eqnarray*}
where $\widetilde{\phi}$ solves $-\Delta
\widetilde{\phi}=\left(\sum_{i=1}^K \lambda_i z_{P_i}\right)
z_{\bf P}$. Reasoning as in the proof of Proposition
\ref{soluzioniapprossimanti}, we obtain that $(A)=o(1)$,
$(B)=o(1)$, $(C)=o(1)$, $(D)=o(1)$. Moreover
\begin{eqnarray*}
(E)\leq C(\lambda_i)\int_{\mathbb R^3}\left[\left|z_{\bf P}\right|^{p+1}-\sum_{i=1}^K z_{P_i}^{p+1}\right]\, dx.
\end{eqnarray*}
Then $(E)=o(1)$ as $\e\rightarrow 0$ (see Proposition
\ref{soluzioniapprossimanti}). At the end
$$I_{\e}''(z_{\bf P})[u, u]=\sum_{i=1}^K\lambda_i^2I''_0(z_{P_i})[z_{P_i}, z_{P_i}]+o(1).$$
Therefore, using Lemma \ref{nondeg} we have, for $\e$ small, that
$$I_{\e}''(z_{\bf P})[u, u] \leq  (1-p)\sum_{i=1}^K\lambda_i^2\|z_{P_i}\|^2 < -C_1<0.$$ So $I_{\e}''(z_{\bf P})$ is negative definite on $A$. We now prove that $I_{\e}''(z_{\bf P})$ is positive definite on $B$. \\
Choose an arbitrary $u\in B$. For simplicity, assume that
$\|u\|=1$. We denote by $\hat{\phi}$ the solution of $-\Delta
\hat{\phi}=z_{\bf P}u$. Since $z_{\bf P}$ and $u$ are bounded, it
is easy to see that, for $\e$ small enough, $$\e^2\int_{\mathbb
R^3}\left[\phi_{z_{\bf P}}u^2+2\hat{\phi}z_{\bf
P}u\right]dx=\int_{\mathbb R^3}\left[V(\e x)-1\right]u^2dx=o(1).$$
Then
\begin{eqnarray*}
I_{\e}''(z_{\bf P})[u, u]&=& \int_{\mathbb R^3}\left[|\nabla u|^2+V(\e x)u^2+\e^2\phi_{z_{\bf P}}u^2+2\e^2\hat{\phi}z_{{\bf P}}u-pz_{\bf P}^{p-1}u^2\right]\, dx\\\\
&=&\int_{\mathbb R^3}\left[|\nabla u|^2+u^2-pz_{\bf P}^{p-1}u^2\right]\, dx+o(1).
\end{eqnarray*}
As done in Proposition \ref{soluzioniapprossimanti} it can be proved that
$$\int_{\mathbb R^3}z_{\bf P}^{p-1}u^2dx= \int_{\mathbb R^3}\sum_{i=1}^Kz_{P_i}^{p-1}u^2dx+o(1).$$ Hence
\begin{eqnarray}\label{Isecondozp}
I_{\e}''(z_{\bf P})[u, u]&=& \int_{\mathbb R^3}\left[|\nabla
u|^2+u^2-p\sum_{i=1}^Kz_{P_i}^{p-1}u^2\right]\, dx+o(1).
\end{eqnarray}
We need to estimate the integral in (\ref{Isecondozp}). In order to do this, we use the following technical result:\\\\
{\it Claim:} for $\e$ small there exists $R\in \left(\e^{\frac{\theta}{2}},
\frac{1}{2}\e^{\theta}\right)$, with
$\theta=\frac{2-\alpha}{\alpha+1}+\delta<0$, such that
\begin{equation}\label{claim}
\sum_{i=1}^K\int_{R<|x-P_i|<R+1}[|\nabla u|^2+u^2]\, dx<4\e^{-\theta}.\\\\
\end{equation}
To prove this we remark that from $\|u\|=1$ it follows
\begin{equation*}
\sum_{i=1}^K\sum_{R\in (\e^{\frac{\theta}{2}},
\frac{1}{2}\e^{\theta})}\int_{R<|x-P_i|<R+1}[|\nabla u|^2+u^2]\,
dx\leq1\qquad R\in \mathbb N.
\end{equation*}
Since, for $\e$ small, the above sum has more than $\frac{\e^{\theta}}{4}$ summands, then, it is always possible to choose $R\in \mathbb N$, $R\in \left(\e^{\frac{\theta}{2}}, \frac{1}{2}\e^{\theta}\right)$ such that the claim holds. \\\\
Let us fix $R$ such that (\ref{claim}) is satisfied and define the smooth cut-off functions $\chi_i: \mathbb R \rightarrow [0, 1]$, $i=1,..., K$ by setting
\begin{equation*}
\chi_i(x):=\left\{
  \begin{array}{lr}
    1 \,\ \,\ \,\ \,\, \qquad\qquad\qquad \,\, |x-P_i|<R
    \\\\
    0\,\ \,\ \,\ \,\, \qquad\qquad\qquad \,\, |x-P_i|>R+1
    \\\\
    |\nabla\chi_i(x)|\leq 2 \qquad\quad \forall\,\, x\in \mathbb R^3.
  \end{array}
 \right.
\end{equation*}
Define also $\chi_0(x)=1-\displaystyle\sum_{i=1}^K \chi_i(x)$. Then we can decompose $u=\displaystyle\sum_{i=0}^K u_i$ where $u_i=u\chi_i$. From (\ref{claim}) it follows that for $i\neq j$ $(u_i, u_j)=o(1)$. Thus
$$1=\|u\|^2=\displaystyle\sum_{i=0}^K \|u_i\|^2+o(1).$$ Using again (\ref{claim}), we obtain that $(z_{P_i}, u_j)=o(1)$ for $i>0$, $i\neq j$. Since $u\in B$ we have $(z_{P_i}, u)=o(1)$. Then for $i=1,..., K$ $$(z_{P_i}, u)=\sum_{j=0}^K(z_{P_i},
u_j)=(z_{P_i}, u_i)+o(1).$$ Hence $(z_{P_i}, u_i)=o(1)$. Finally, for $i=1,...,K$, since $u\bot
\dot{z}_{i, j}$, reasoning as above, we find also $(\dot{z}_{i,
j}, u_l)=o(1)$ for all $i, j, l$.\\
By using the above properties and Lemma \ref{nondeg} we obtain
\begin{eqnarray*}
I_{\e}''(z_{\bf P})[u, u]&=& \int_{\mathbb R^3}\left[|\nabla u|^2+u^2-p\sum_{i=1}^Kz_{P_i}^{p-1}u^2\right]\, dx+o(1)\\\\
&=&\sum_{i=1}^K \int_{\mathbb R^3}\left[|\nabla u_i|^2+u_i^2-pz_{P_i}^{p-1}u_i^2\right]\, dx+\|u_0\|^2+o(1)\\\\
&\geq & C\sum_{i=1}^{K} \|u_i\|^2+\|u_0\|^2+o(1)\\\\
&\geq & C_2\left(\sum_{i=0}^{K}\|u_i\|^2\right)+o(1)\\\\
&\geq & C_2 >0.
\end{eqnarray*}
\end{Proof}
With this estimates in hand we can now solve the auxiliary
equation. Consider $z= z_{\bf P} \in \mathcal{Z}$ fixed, and
define
$$B_{\e}=\left\{ u\in W : \|u\|\leq 2
\bar{C} \|I'_{\e}(z)\|\right\},$$ where $\bar{C}$ is the positive
constant given by Proposition \ref{prop2}. So, the solutions of
the auxiliary equations are fixed points of the map
$S_{\e}:W\rightarrow W$
$$S_{\e}(w)=w-[PI_{\e}''(z)]^{-1}[PI'_{\e}(z+w)].$$
It is easy to check that $\| S_{\e}(0)\| \leq \bar{C}
\|I'_{\e}(z)\|$. We now compute the derivative of $S_{\e}$:
$$S_{\e}'(w)[v]= v- [PI_{\e}''(z)]^{-1}PI_{\e}''(z+w)[v] =
[PI_{\e}''(z)]^{-1}\left (PI_{\e}''(z)-PI_{\e}''(z+w) \right
)[v].$$
Now observe that $I''_{\e}$ is uniformly continuous in bounded
sets, so
$$\|PI_{\e}''(z+w)-PI_{\e}''(z)\| \rightarrow 0 \qquad
(\e\rightarrow 0) $$ uniformly in $z\in \mathcal{Z}$ and $w\in
B_{\e}$ (recall Proposition \ref{soluzioniapprossimanti}).

\medskip This implies that $\| S_{\e}'(w)\|= o(1) $ for any $w \in
B_{\e}$. Therefore, $S_{\e}$ is a contraction and, by using the
mean value theorem, $S_{\e}(B_{\e}) \subset B_{\e}$. We make use
of the Banach contraction theorem to find a unique fixed point
$w=w_{\e,z}\in B_{\e}$ of $S_{\e}$. Moreover one has
\begin{equation}\label{stimawminimo}
\|w_{\e, z}\| \leq 2\bar{C}\|I'_{\e}(z)\| \leq C \e^2
\end{equation}

\medskip

\section{The reduced functional} \label{sezioneridotto}
In this section we will find a solution for the bifurcation
equation among the set of solutions of the auxiliary equation,
which is:
$$\bar{\mathcal{Z}}=\left\{ z+w_{\e,z} : z \in \mathcal{Z},\ w_{\e, z} \mbox{
solves } \eqref{lya} \mbox{(a)}, \mbox{ and satisfies }
\eqref{stimawminimo} \right \}.$$

By the Implicit Function Theorem it is easy to check that
$\bar{\mathcal{Z}}$ is a $C^1$ manifold. Moreover, it is
well-known (see \cite{AmbrMalc}, for example) that
$\bar{\mathcal{Z}}$ is a natural constraint for $I_{\e}$ for $\e$
small. In other words, critical points of
$I_{\e}|_{\bar{\mathcal{Z}}}$ are solutions of the bifurcation
equation \eqref{lya} (b), and hence solutions of
\eqref{problemafinale}.

So, let us define the reduced functional as the restriction of the functional $I_{\e}$ to the natural constraint $\bar{\mathcal{Z}}$, namely $\Phi_{\e}: \Lambda_{\epsilon} \to \R$,
$\Phi_{\e}({\bf P})=I_{\e}(z_{\bf P}+w_{\e, z_{\bf P}})$, and we look for critical points of
$\Phi_{\e}$. Using the information on $\|w_{\e, z_{\bf P}}\|$, we will be able to
find an expansion of $\Phi_{\e}({\bf P})$.\\
First of all, since $I_{\e}''$ maps bounded sets onto bounded
sets, we have

\begin{equation*}\label{ridottoprimo0}
\Phi_{\e}({\bf P})=I_{\e}(z_{\bf P})+I'_{\e}(z_{\bf
P})[w_{\e, z_{\bf P}}]+O(\|w_{\e, {\bf P}}\|^2).
\end{equation*}
Using Proposition \ref{soluzioniapprossimanti} and (\ref{stimawminimo}) we deduce

\begin{equation}\label{ridottoprimo}
\Phi_{\e}({\bf P})=I_{\e}(z_{\bf P})+O(\e^4).
\end{equation}
So we have to compute $I_{\e}(z_{\bf P})$. Preliminary lemmas are
in order.

\begin{lemma}\label{lemmaperespansionephi}
For $\beta =1, 2$ and $F: \mathbb R^3 \rightarrow \mathbb R$ such
that $(1+|y|^{\beta+1})F\in L^1 \cap L^{\infty}$ set
$$\Psi_{\beta}[F](x)=\int_{\mathbb
R^3}\frac{1}{|x-y|^{\beta}}F(y)\, dy.$$ Then there exist two
positive constants $C=C(\beta, F)$ and $C'=C'(\beta, F)$ such that
\begin{equation}\label{espphi}
\left|\Psi_{\beta}[F]-\frac{C}{|x|^{\beta}}\right|\leq
\frac{C'}{|x|^{\beta+1}}, \qquad \forall\,\, x \neq 0.
\end{equation}
\end{lemma}

For a proof see \cite{DaprileWei}.\\ Now, thanks to the the
exponential decay of $U$ the following estimate holds (see Lemma
2.1 of \cite{KangWei}):

\begin{lemma}\label{perpezzimisti}
For $\e$ sufficiently small and ${\bf P} \in \Lambda_{\e}$, we
have
$$\int_{\mathbb R^3}z_{P_i}^pz_{P_j}\, dx=(\eta+o(1))e^{-|P_1-P_2|}$$
where $$\eta=\int_{\mathbb R^3}U^p(x)e^{-x_1}\, dx>0.$$
\end{lemma}
We are now in position to find an expansion of $I_{\e}(z_{\bf P})$.

\begin{proposition}\label{Iepsilonzp}
For any ${\bf P}=(P_1, ..., P_K) \in \Lambda_{\e}$ and $\e>0$
sufficiently small we have

\begin{equation}\label{ridottosecondo}
I_{\e}(z_{\bf P})= C_0 +\e^2 C_1+C_2\sum_{i=1}^K V(\e P_i)+C_3
\e^2 \sum_{i\neq j}\frac{1}{|P_i-P_j|}
+o(\e^{\frac{3\alpha}{\alpha+1}-\delta})
\end{equation}

where $$C_0=K\cdot\left(\frac{1}{2}\int_{\mathbb R^3}|\nabla
U|^2\, dx-\frac{1}{p+1}\int_{\mathbb R^3}U^{p+1}\, dx\right),$$
$$C_1=\frac{K}{4}\int_{\mathbb R^3}\frac{U^2(x)U^2(y)}{|x-y|}\,
dx\, dy,\qquad C_2=\frac{1}{2}\int U^2\, dx,$$  and $C_3$ is a positive constant given by Lemma \ref{lemmaperespansionephi},
which depends only on $U$.

\end{proposition}
\begin{Proof}
We compute
\begin{eqnarray*}
I_{\e}(z_{\bf P})&=&\sum_{i=1}^K I_{\e}(z_{P_i})+\sum_{i\neq j}\int_{\mathbb R^3}\left[\nabla z_{P_i}\nabla z_{P_j}+V(\e x)z_{P_i}z_{P_j}\right]dx+\frac{\e^2}{4}\sum_{i\neq j}\int_{\mathbb R^3}\phi_{z_{P_i}}z_{P_j}^2dx\\
&&+\frac{\e^2}{2}\sum_{l, i\neq j}\int_{\mathbb R^3}\phi_{i, j}z_{l}^2dx+\frac{\e^2}{4}\sum_{i\neq j}\int_{\mathbb R^3}\phi_{z_{\bf P}}z_{P_i}z_{P_j}dx\\
&&-\frac{1}{p+1}\int_{\mathbb R^3}\left[\left|z_{\bf P}\right|^{p+1}-\sum_{i=1}^K|z_{P_i}|^{p+1}\right]dx\\
&=&\sum_{i=1}^K I_{\e}(z_{P_i})+\sum_{i\neq j}\int_{\mathbb R^3}z_{P_i}^pz_{P_j}dx+\frac{\e^2}{4}\sum_{i\neq j}\int_{\mathbb R^3}\phi_{z_{P_i}}z_{P_j}^2 dx\\
&&+\frac{\e^2}{2}\sum_{l, i\neq j}\int_{\mathbb R^3}\phi_{i, j}z_{l}^2dx+\frac{\e^2}{4}\sum_{i\neq j}\int_{\mathbb R^3}\phi_{z_{\bf P}}z_{P_i}z_{P_j}dx\\
&&-\frac{1}{p+1}\int_{\mathbb R^3}\left[\left|z_{\bf P}\right|^{p+1}-\sum_{i=1}^K|z_{P_i}|^{p+1}\right]dx+o(\e^{\frac{3\alpha}{\alpha+1}-\delta}).\\
\end{eqnarray*}
Here $\phi_{i, j}$ are the solutions of $-\Delta \phi_{i,
j}=z_{P_i}z_{P_j}$, $i\neq j$. Let us evaluate separately the
various terms.

{\bf Claim:} There holds:
\begin{equation}\label{sommaizpi}
 I_{\e}(z_{P_i})=\widetilde{C}_0+\e^2 \widetilde{C_1}+C_2V(\e P_i)+o(\e^{\frac{3\alpha}{\alpha+1}-\delta})
\end{equation}
where $$\widetilde{C}_0=\frac{1}{2}\int_{\mathbb R^3}|\nabla
U|^2dx-\frac{1}{p+1}\int_{\mathbb R^3}|U|^{p+1}dx,\quad
\widetilde{C_1}=\frac{1}{4}\int_{\mathbb R^3}\phi_{U}U^{2}dx,\quad
C_2=\frac{1}{2}\int_{\mathbb R^3}U^2dx.$$ \\

It suffices to estimate:
$$ \int_{\R^3} [V(\e x)- V(\e P_i)] U^2(x-P_i) \, dx. $$

First, we split this integral expression in two terms

\begin{eqnarray*} \int_{\R^3}[V(\e x)-V(\e P_i)]z_{P_i}^2\, dx = &
\dis \int_{|x-P_i|>\e^{-\tau}}[V(\e x)-V(\e P_i)]z_{P_i}^2\, dx + \\
& \dis \int_{|x-P_i|<\e^{-\tau}}[V(\e x)-V(\e P_i)]z_{P_i}^2\, dx,
\end{eqnarray*} for some positive constant $\tau$ to be determined.
Since $V$ is bounded in $L^{\infty}$, we use the change $y=x-P_i$,
and the exponential decay of $U$ to conclude
$$\int_{|x-P_i|>\e^{-\tau}}[V(\e x)-V(\e
P_i)]z_{P_i}^2\, dx \leq C \int_{|y|>\e^{-\tau}}U^2(y)\, dy  =
o(\e^M)$$ for any positive $M$.

We use a Taylor expansion:

\begin{equation} \label{taylor} \begin{array}{c} \left |
\dis \int_{|x-P_i|< \e^{- \tau}} \left [V(\e x)-V(\e P_i) - \e
\nabla V(\e P_i) \cdot
(x-P_i) \, \right ] z_{P_i}^2 \right | \leq  \\
\dis \frac{\e^2}{2} \max \{ \|D^2V(\xi)\|: \ |\xi - \e P_i| <
\e^{1-\tau} \} \, \dis \intr |x-P_i|^2 z_{P_i}^2 \, dx.
\end{array}
\end{equation}

By using the radial symmetry of $U$,

$$ \int_{|x-P_i|< \e^{- \tau}} \nabla V(\e P_i) \cdot
(x-P_i) U(x-P_i)^2 =0. $$

So, it suffices to estimate $\|D^2V(\xi)\|$ for $\ |\xi - \e P_i|
< \e^{1-\tau}$. First, observe that if $\tau<1$ and $\e$ is small
enough, $\xi \in \mathcal{U}$.

Moreover, by the definition of $\Lambda_{\e}$, $V(\e P_i) = 1 +
|g(\e P_i)|^{\alpha}\leq 1+\e^{\frac{3\alpha}{\alpha+1}-\delta}$.
From this and \eqref{second} we have that $$|V_{x_i \, x_j}(\e
P_i)| \leq C |g(x)|^{\alpha -2} \leq C \e^{
\frac{\alpha-2}{\alpha} \left ( \frac{3\alpha}{\alpha+1}-\delta
\right ) }.$$

On the other hand, since $V \in C^{2,\gamma}$ (recall, $\gamma=
\min\{1, \alpha-2\}$):

$$ \left | V_{x_i \, x_j}(\xi)-V_{x_i \, x_j}(\e P_i) \right | \leq C \e^{\gamma (1-\tau)}.$$

Therefore,
$$ | V_{x_i \, x_j}(\xi)| \leq C \e^{\min \{\frac{\alpha-2}{\alpha}\left (
\frac{3\alpha}{\alpha+1} - \delta \right ),\ \gamma (1-\tau)\} }.
$$

By direct computation,  $2+ \frac{\alpha-2}{\alpha}\left (
\frac{3\alpha}{\alpha+1} - \delta \right ) >
\frac{3\alpha}{\alpha+1} - \delta$. Moreover, $2+1=3 >
\frac{3\alpha}{\alpha+1}$ and $2+ \alpha -2= \alpha >
\frac{3\alpha}{\alpha+1}$. Then, we can choose $\tau>0$ small
enough such that $2+ \gamma (1-\tau)>
\frac{3\alpha}{\alpha+1}-\delta$. This concludes the proof of the
claim.

\medskip

We now continue the estimates of the remaining terms. From Lemma
\ref{perpezzimisti}
\begin{equation}\label{zpipzpj}
 \int_{\mathbb
 R^3}z_{P_i}^p z_{P_j}dx=(\eta+o(1))e^{-|P_i-P_j|}=o(\e^M)
\end{equation}
for any $M>0$. Now, by using the notations of Lemma
\ref{lemmaperespansionephi}, we have
$\phi_{z_{P_i}}=\frac{1}{4\pi}\Psi_1[U^2](x-P_i)$. If $i\neq j$,
by (\ref{espphi})
\begin{eqnarray*}
\int_{\mathbb R^3}\phi_{z_{P_i}}z_{P_j}^2\, dx&=&\frac{1}{4\pi}\int_{\mathbb R^3}\Psi_1[U^2](x-P_i)U^2(x-P_j)\, dx\\\\
&=& C_3\int_{\mathbb R^3}\frac{1}{|x-P_i|}U^2(x-P_j)\, dx+O(1)\int_{\mathbb R^3}\frac{1}{|y+P_j-P_i|^2}U^2(y)\, dx\\\\
&=& C_3\Psi_1(U^2)|P_i-P_j|+O(1)|P_i-P_j|^{-2}.
\end{eqnarray*}
From the definition of $\Lambda_{\e}$ and since $\alpha >2$, $|P_i -
P_j|^{-2} = o(\e^{\frac{3 \alpha}{\alpha +1}-\delta})$.
Furthermore
$$\int_{\mathbb R^3}\phi_{z_{\bf P}}z_{P_i}z_{P_j}\, dx\leq C\int_{\mathbb R^3} z_{P_i}z_{P_j}dx = o(\e^M) \  (i\neq j).$$
and, consequently,
\begin{equation*}
\int_{\mathbb R^3}\phi_{i, j}z_{P_i}^2\, dx=-\intr \phi_{i,j}
\Delta \phi_{z_{P_i}} = -\intr \Delta \phi_{i,j} \,
\phi_{z_{P_i}}= \int_{\mathbb R^3}\phi_{z_{P_i}}z_{P_i}z_{P_j}\,
dx=o(\e^M)
\end{equation*}
for any $M>0$. Since ${\bf P}\in \Lambda_{\e}$, we have that for
$i\neq j$
\begin{equation*}
 \frac{\e^2}{4}\int_{\mathbb R^3}\phi_{z_{P_i}}z_{P_j}^2\, dx =C_3\frac{\e^2}{|P_i-P_j|}+\e^2 O(|P_i-P_j|^{-2})=
 C_3\frac{\e^2}{|P_i-P_j|}+ o(\e^{\frac{3 \alpha}{\alpha
 +1}-\delta}).
\end{equation*}

Finally, arguing as in Proposition \ref{soluzioniapprossimanti}, we obtain
\begin{equation}\label{ultimo}
\int_{\mathbb R^3}\left[\left|z_{\bf P}\right|^{p+1}-\sum_{i=1}^K
|z_{P_i}|^{p+1}\right]\, dx=o(\e^{M})
\end{equation} for any $M>0$.

All previous estimates imply the expansion
\eqref{ridottosecondo}.

\end{Proof}

From (\ref{ridottoprimo}) and (\ref{ridottosecondo}) we have the
following expansion for the reduced functional
\begin{equation}\label{ridottofinale}
\Phi_{\e}({\bf P})=C_0 +\e^2 C_1+C_2\sum_{i=1}^K V(\e P_i)+C_3
\e^2 \sum_{i\neq
j}\frac{1}{|P_i-P_j|}+o(\e^{\frac{3\alpha}{\alpha+1}-\delta}).
\end{equation}

\begin{proposition}\label{risoluzione}
For $\e$ sufficiently small, the following minimization pro-\ blem
\begin{equation}\label{minimizzazione}
\min\left\{ \Phi_{\e}({\bf P})\,\ : \,\ {\bf P} \in \Lambda_{\e}
\right\}
\end{equation}
has a solution ${\bf P}_{\e}\in \Lambda_{\e}$.
\end{proposition}
\begin{Proof}
Since $\Phi_{\e}(\bf P)$ is continuous in $\bf P$ in a compact set,
the minimization problem has a solution. Let $\Phi_{\e}(\bf
P^{\e})$ be the minimum of $\Phi_{\e}$ where ${\bf P}^{\e}$ is in the closure of the set $\Lambda_{\e}$. We prove by energy
comparison that $\bf P^{\e}$ is not on the boundary of
$\Lambda_{\e}$. In order to do this, first we obtain an upper bound
for $\Phi_{\e}(\bf P^{\e})$. Let us choose
$$P_{j}^0=\e^{\frac{2-\alpha}{\alpha+1}}X_j$$ where $X_j,
\,\ j=1,..., K$ are the $K$ vortices of $K-$polygon centered at
$0$ with $|X_i-X_j|=1$, $i\neq j$.\\ Then for $\e$ small it is
clear that $\e P_j^0 \in \mathcal{U}$. Moreover
$$|P_i^0-P_j^0|=\e^{\frac{2-\alpha}{\alpha+1}}|X_i-X_j|$$ and
$$V(\e P_j^0)\leq 1+C |\e
P_j^0|^{\alpha} \leq 1+ C \e^{\frac{3\alpha}{\alpha+1}}.
$$ Therefore, ${\bf
P^0}=(P_1^0,..., P_K^0) \in \Lambda_{\e}$. Hence by
(\ref{ridottofinale}) we obtain
\begin{equation}\label{upperbound}
\Phi_{\e}({\bf P^{\e}})=\min_{{\bf P}\in
\Lambda_{\e}}\Phi_{\e}({\bf P})\leq \Phi_{\e}({\bf P^0})\leq C_0 +
\e^2 C_1 + K C_2+  C_3 \e^{\frac{3\alpha}{\alpha+1}}.
\end{equation}
If now ${\bf P^{\e}}$ is such that
$|P_i^{\e}-P_j^{\e}|=\e^{\frac{2-\alpha}{\alpha+1}+\delta}$ for
some $i\neq j$, then
\begin{equation}\label{primoassurdo}
\Phi_{\e}({\bf P^{\e}})\geq C_0 + \e^2 C_1 + K C_2 + C_3
\e^{\frac{3\alpha}{\alpha+1}-\delta}.
\end{equation}
If, instead, ${\bf P^{\e}}$ is such that $V(\e
P_i^{\e})=1+\e^{\frac{3\alpha}{\alpha+1}-\delta}$ for some $i$,
then
\begin{equation}\label{secondoassurdo}
\Phi_{\e}({\bf P^{\e}})\geq C_0 + \e^2 C_1 +K C_2 + C_2
\e^{\frac{3\alpha}{\alpha+1}-\delta}.
\end{equation}
But both (\ref{primoassurdo}) and (\ref{secondoassurdo}) are in contradiction with (\ref{upperbound}).\\

We remark that we have not considered the case $\e {\bf P}^{\e}\in
\partial \mathcal{U}$, because this would be in contradiction with
$V(\e P_{j}^{\e}) \leq 1+\e^{\frac{3\alpha}{\alpha+1}-\delta}$ for
$\e$ small.
\end{Proof}

\end{document}